\documentclass[a4paper]{article}
\usepackage{amsmath, amssymb}
\usepackage{stmaryrd}
\usepackage{amsthm,cite}
\usepackage{amsfonts,mathtools}
\usepackage{graphicx,subfigure}
\usepackage{epsfig,color,pifont,enumerate}
\usepackage{setspace}
\graphicspath{{img/}}

\usepackage{geometry}
\renewcommand{\emptyset}{\O}
\usepackage[mathscr]{eucal}
\renewcommand{\phi}{\varphi}
\renewcommand{\epsilon}{\varepsilon}
\renewcommand{\kappa}{\varkappa}
\newcommand{\sgn}{\text{\bf sgn}\,}

\newcommand{\br}{\mathbb{R}}
 \newcommand{\bldr}{\boldsymbol{r}}

\newtheorem{theorem}{Theorem}[section]

\newtheorem{define}{Definition}[section]

\newtheorem{assumption}{Assumption}[section]
\newtheorem{lemma}{Lemma}[section]

\newtheorem{proposition}{Proposition}[section]

\newcommand{\spr}[2]{\left\langle #1; #2 \right\rangle}

\newcommand{\ov}[1]{\overline{#1}}

\newcommand{\ve}{\varepsilon}

\newcommand{\pf}{{\bf Proof:}\;}
\newcommand{\epf}{$\qquad \Box$}

\newcommand{\imm}{\vec{\boldsymbol{\imath}}}
\newcommand{\col}{\mathbf{col}}
\newcommand{\sprl}[2]{\langle #1; #2 \rangle}
\newcommand{\dsafe}{d_{\text{safe}}}
\newcommand{\ts}{\mathscr{T}}
\newcommand{\ff}[3]{\mathbf{II}_{#1}\!\left[ #2; #3 \right]}
\newcommand{\sprs}[2]{\langle #1; #2 \rangle}
\newcommand{\bro}{\boldsymbol{\varrho}}
\title{Proofs of the Technical Results Justifying an Algorithm of Reactive 3D Navigation for a Surface Scan by a  Nonholonomic Mobile Robot\footnote{This work was supported by the Australian Research Council, and by the Russian Science Foundation 14-21-00041 and the Saint Petersburg State University.}}
\date{}
\author{Alexey S. Matveev\footnotemark[2], Kirill S. Ovchinnikov\footnotemark[2], and Andrey V.Savkin\footnotemark[3]}
\begin{document}
\maketitle
\renewcommand{\thefootnote}{\fnsymbol{footnote}}
\footnotetext[2]{Department of Mathematics and Mechanics, Saint Petersburg State University, St. Petersburg, Russia}
\footnotetext[3]{School of Electrical Engineering and Telecommunications, University of New South Wales, Sydney,  Australia}
\section{Introduction}
The paper addresses autonomous navigation of a mobile robot in missions such as aerial and underwater photography, environmental monitoring, remote sensing, border patrol, environmental monitoring, imaging and inspecting various 3D structures.
For many of these and other missions, it is needed to drive a
robot in a 3D space so that it advances to a certain 2D surface, then
moves close to it and performs its scanning, i.e., sequentially visits tight vicinities of all its points.
The seminal paper \cite{gage:92} classifies such behavior as {\it sweep coverage} of the surface.
Apart from the afore-mentioned examples, navigation of mobile robots for sweep coverage
may be of interest for environmental
studies, geological exploration, underwater archeology \cite{BiFoSiCa10}, seabed exploration \cite{HeTiLu92}, rescue and search operations, fine-scale imaging and inspection of underwater, air, or space structures, e.g., coral reefs or ship hulls
\cite{JoRoPiWi10,HoEuKiEn12}, spray painting \cite{AtGrCoCh05}, oil and gas pipeline inspection, infrastructure surveillance, damage assessment, recovery, and protection of 3D objects, localization, reduction, and stopping of underwater oil spills \cite{SaChXiJaMa15}, 3D building reconstruction, and aerial city surveillance, to name just a few.
\par
Sweep coverage of planar objects in a two-dimensional workspace, also known as area coverage, enjoys a substantial literature; we refer the reader to
\cite{Choset01,GalCar13,HsLiYa14,SaChXiJaMa15} for comprehensive surveys.
Most of these works assume that both the area to be covered and robot reside in a common plane. However, this fails to be true for many applications, e.g., for inspecting the surfaces of chimney-stalks or gas storage tanks. Meanwhile, only a few works address special challenges arisen in coverage of essentially non-planar surfaces when moving in a three-dimensional workspace.
\par
Vertically projectively planar surfaces \cite{HeTiLu92} can be handled via direct application of planar terrain-covering algorithms (e.g., lawnmover-like patterns) to the robot's vertical projection onto a horizontal plane with either following the elevation profile of the terrain within a certain range of depths \cite{HeTiLu92} or successively moving in evenly spaced parallel planes \cite{LeChLe09}.
In \cite{AtChRiAc01}, it is advocated to follow the horizontal slice of a closed orientable surface with a pre-specified margin and to vertically lift, from time to time, the slicing plane by a given step. Critical points are detected to partition a topologically complex surface into simple parts to be individually scanned. The problem of coverage path planning directly on the surface is addressed in \cite{AtGrCoCh05,AtGrCoCh09} based on a cellular decomposition of the surface into geometrically simple patches (extruded surfaces homeomorphic to a disk), which are subjected to relatively independent scans. This opens the door for primarily ``planar'' treatment of any patch, whereas performance issues may cause the need to allow for surface contortion. A relevant problem of finding an optimal family of curve segments on a surface is considered in \cite{KiSa03}.
Paper \cite{GalCar13} recommends segmentation of the seabed into low-relief and high-slope subregions based on a bathymetric map;
the former are handled via classical lawnmower survey, the latter are handled by following constant-depth horizontal bathymetric contours at iteratively increasing depths.
A somewhat similar approach to coverage path planning for arable farming is employed in \cite{JiTa11}, where
the terrain is partitioned into subregions based mainly on the slope steepness.
\par
Cluttered workspaces and surfaces with occluded parts visible only from special viewpoints are typically treated via consecutive solution
of two NP-hard problems \cite{DaKa00,EngHo12,PaKuPa13}.
First, a kind of the standard art gallery problem is solved in order to determine a set of viewpoints that establish full visual coverage of the surface. Second, a variant of the traveling salesman problem is solved to acquire a shortest path visiting all these points.
The accompanying hard computational burden can be somewhat eased by using a random sampling machinery \cite{DaKa00,EngHo12}.
\par
Thus, the entire previous research is focused on geometrical path planning, implementation of following the planned path
was out of discussion. However, this implementation may be a real problem in the face of practical limitations of the robot; examples are concerned with, but are not limited to, speedy, yet safe, following along a highly contorted surface or reliable estimation of the surface's slope based on noisy data. Moreover, kinematic, dynamic, or sensing constraints may cause infeasibility of
the planned maneuvers or operations, like a sharp vertical turn to escape to a higher plane or detecting the end of a full turn around a horizontal slice of the surface. When solving purely navigation tasks, most of the previous works either hinge on a known map of the environment, which may be unavailable in practice, or involve computationally expensive sensory data processing, close to building such a map.
\par
On the other hand, there is a growing demand to perform coverage operations faster, yet with improved safety and efficiency, e.g., closer to the surface, and by means of simpler and cheaper hardware. This calls for computationally inexpensive algorithms
that allow a mobile robot to quickly make and implement navigation decisions and ensure safety, while respecting all perceptional and computational limitations and motion constraints that may apply.
\par
This paper is aimed at demonstration that even in the face of challenges from nonholonomic constrains, under-actuation, and finite control range, safe sweep coverage of a 2D surface embedded in a 3D workspace can be achieved at the lowest control level via a computationally primitive reflex-like rules. They straightforwardly transform the current sensor readings into a current control input and do not use complex or dubious operations like even short-term motion planning, building a map of even a tight vicinity of the scene or surface, or depositing marks in the environment.
Moreover, the above challenges are enhanced by poor perceptual capacity of the robot that is confined to access only to a certain direction in 3D, to the robot's coordinate along it, and to the distance from the scanned surface in the perpendicular plane. Finally, no a priori knowledge about the surface is available and the memory size of the robot is drastically small so that it is impossible to memorize, in any form, online data about the surface in even ``next to nothing'' portions.
\par
To the above end, we address control at the kinematic level of a robot described by the standard model of the 3D nonholonomic unicycle:
The robot travels at a constant speed in 3D and is steered by the acceleration vector, which is bounded in magnitude and is constantly normal to the velocity. This classic model applies to many mechanical systems, e.g., constant-speed vehicles moving in the surge direction and controlled by bounded yaw and pitch rates, like fixed-wing aircrafts or torpedo-like underwater vehicles. Moreover, this model is relevant for
any systems whose acceleration can be freely manipulated within a disk normal to the velocity vector, like for many rotorcraft systems \cite{CaChLe11}. The objective is to perform sweep coverage of the surface via spiraling around it with a drift in the assigned direction.
\par
In this paper, we deal with compact surfaces that are not normal to this direction within the scan range.
At the same time, the respective findings may be used as primitives when handling more complex surfaces via their segmentation.
To start with, we disclose a trade-off between the surface's contortion and the limited turning capacity of the robot necessary for the requested  sweep coverage to be feasible. After this, we assume that these conditions are satisfied in a slightly enhanced form, propose a novel navigation
law, and show, via a mathematically rigorous nonlocal convergence result, that the control objective is achieved by this law provided that it is properly tuned. Thus this law solves the problem almost whenever it is solvable. Explicit recommendations on tuning the controller are also offered.
\par
This paper develops some ideas from \cite{MTS11,MTS11ronly,MaTeSa11,MaTeSa12,MaChaSa12a,MaHoSaRAS13,MaHoSa3d14,HoMaSa15rob,MaSaHoWa16b}.
The form of the 3D kinematic model is borrowed from \cite{MaHoSa3d14}, where its detailed discussion, including relations with a more conventional form, is offered. The other papers handle 2D environments
and are threaded by a general regulation paradigm stemming from the biologically inspired equiangular navigation law \cite{TS10}.
However, the 3D context causes serious extra challenges and makes the findings of \cite{MTS11,MTS11ronly,MaTeSa11,MaTeSa12,MaChaSa12a,MaHoSaRAS13,MaSaHoWa16b} mostly inapplicable. For example, though the targeted ``spiralling'' can be viewed as 2D motion about the surface combined with a perpendicular drift, the properties of 2D motion are far from key those of a Dubins car treated in the above papers, to say nothing that the respective loops of regulation of these 2D motion and drift are not decoupled, whereas the influence of the latter on the former carries a potential of failure and calls for thorough extra research.
\par
The body of the paper is organized as follows. Section~\ref{sec1}
describes the problem setup, its technical specification is give in Section~\ref{sec.specif1}.
Section~\ref{sec.specif} is devoted to necessary conditions and associated assumptions of theoretical analysis.
The control law and the main results are presented in Sections~\ref{sec.claw}. These results are illustrated in Section~\ref{sec.exmpl} for special scenarios.
\par
\par
The extended introduction and discussion of the proposed control law are given in the paper submitted by the authors to the IFAC journal {\it Automatica}. This text basically contains the proofs of the technical facts underlying justification of the convergence at performance of the proposed algorithm in that paper, which were not included into it due to the length limitations. To make the current text logically consistent, we reproduce the problem statement and notations.
\section{Description of the System and General Setup of the Problem}
\label{sec1}
A mobile robot travels in the three-dimensional space $\br^3$ with a constant surge speed $v$ and is controlled by bounded pitching $q$ and yawing $r$ rates.
\footnote{The projections  of the angular velocity on the pitch and yaw axes of the
standard robot-fixed Cartesian frame}
There is an unknown domain $D \subset \br^3$ in the space. The objective is to drive the robot to the desired distance $d_0$ to $D$ and to scan and patrol its boundary surface $\partial D$ afterwards, with maintaining this distance $d_0$.
\par
By following \cite{MaHoSa3d14}, we do not come into details of the roll motion (which is of minor importance for the problem at hand) and employ the abridged kinematic model of 3D unicycle:
\begin{equation}
\label{1b}
\dot{\bldr} = v \vec{\boldsymbol{\imath}}, \quad \dot{\vec{\boldsymbol{\imath}}} = \vec{u}, \quad \spr{\vec{u}}{\imm} =0, \quad \|\vec{u}\|  \leq \ov{u}.
\end{equation}
Here $\bldr = \col(x,y,z)$ is the robot's position, $x,y,z$ are its absolute Cartesian coordinates, $\imm$ is the unit vector along its centerline,
$\spr{\cdot}{\cdot}$ and $\| \cdot \|$ are the standard Euclidian inner product and norm, respectively, $\vec{u}$ is the two-dimensional control, and the upper bound $\ov{u}$ on its magnitude is given. In \eqref{1b}, the equation $\spr{\vec{u}}{\imm} =0$ keeps the length of $\imm$ constant, as is required. In effect, \eqref{1b} comes to setting up the speed $v$ and minimal turning radius $R:=v/\ov{u}$ of the robot: It travels with the speed $v$ over space curves whose curvature radii are no less than $R$.
\par
Remark~2.1 in \cite{MaHoSa3d14} discusses replacement of the ``abstract'' control $\vec{u}$ by the conventional pitching $q$ and yawing $r$ rates in various contexts, based on a one-to-one correspondence $\vec{u} \leftrightarrow (q,r)$.
That remark also shows that the model \eqref{1b} is applicable whenever the speed $v$ can be kept constant by a proper control, whereas the acceleration can be simultaneously manipulated within a disk perpendicular to the velocity $\vec{v}$ and centered at the origin; then $\vec{\imath}:= \vec{v}/v$.  This holds for helicopters, submarine-like vehicles, and many other mechanical systems that move not necessarily in the surge direction \cite{CaChLe11,ReGeChFuLe12}. Since our  main results are expressed in terms of $\vec{u}$, they cover this case as well.
 \par
The robot has access to a certain unit vector $\vec{h} \in \br^3$ in the local frame attached to its body.
It also measures its own coordinate $h(\bldr) = \langle \bldr;\vec{h}\rangle + \text{const}$ along the direction given by $\vec{h}$.
For the convenience's sake and to indicate a typical, though not mandatory, scenario, $\vec{h}$ and $h$ are termed the {\it vertical} vector and {\it altitude}, respectively.
The ``horizontal'' distance $d(\bldr):=\min_{\bldr^\prime \in D \cap H_{\bldr}} \|\bldr - \bldr^\prime\|$ to the domain $D$ is also accessible.
Here $H_{\bldr} := \{\bldr^\prime: \sprl{\bldr^\prime-\bldr}{\vec{h}} =0 \}$ is the horizontal plane passing through the current location $\bldr$ of the robot, and $\inf_{\bldr^\prime \in \emptyset}d:=\infty$. No other measurements are available. A safety margin $d \geq \dsafe > 0$ should be always respected.
 \par
To specify the problem, we note that in our theoretical analysis, $D$ is viewed as a compact $3$-dimensional smooth manifold with boundary. Its boundary surface $\partial D$ is a $2$-dimensional smooth manifold in $\br^3$ with finitely many connected components.
One of them is the ``outer'' boundary of $D$, the others are attached to the ``holes'' in $D$. The objective is to scan a certain component $B_\dagger$, e.g., the ``outer'' one.
The scan should be done within a given range of altitudes $h \in [h_-,h_+]$, where $h_- < h_+$. In the horizontal flat layer $L$ defined by this range, the surface $B_\dagger$ may in turns, look like a set of several isolated pieces; the scan objective is focused on only one of them $B$. So the directive to maintain a distance of $d_0$ in fact addresses the distance to $B$.
\par
The scan itself is meant as a robot's motion for which the point $b_{\bldr}$ of $B$ horizontally nearest to the robot sweeps over the entire $B$ so that any point $b$ of $B$ appears to be in a close proximity of the path traced by $b_{\bldr}$. Then $b$ is at approximately a distance of $d_0 =\|\bldr-b_{\bldr}\|$ from the robot at some time, which distance is viewed as favorable for achieving a certain ``higher level'' objective, like monitoring, recording, processing or inspection of the boundary, or prevention of intrusions through $B$. Finally, repeated scans should be arranged:
a complete surface scan of $B$ should be rerun over and over again.
\par
Now we pass to formal details of the problem setup.
\section{Assumptions Underlying Theoretical Analysis, Necessary Conditions, and Notations}
\label{sec.specif}
This specification proceeds from the following.
\begin{assumption}
\label{ass.smooth}
The domain $D$ is a compact $3$-dimensional $C^3$-smooth manifold with boundary $\partial D$. Within the horizontal layer $L:=\{\bldr: h_- \leq h \leq h_+\}$, the unit inner normal $N_b$ to $\partial D$ is not vertical (co-linear to $\vec{h}$) at any point $b$ from the boundary piece $B_\dagger$ of interest.
\end{assumption}
Here ``inner''$=$``directed inwards $D$''. Due to the second claim, $B_\dagger \cap L$ has finitely many connected components, each vertically spreading from $h_-$ to $h_+$. As was noted, one of them $B$ is to be scanned. Let $D_{\text{out}}$ be the part of the free space in the layer $L$ whose boundary includes $B$.
\par
We would like to organize scan as circumnavigation of $B$ at a distance of $d=d_0$ with a small increment $\Delta h$ in altitude from one round of circumvention to the next one so that the robot spirals around $B$. Circumnavigation may be either outer or inner, depending on the side of $B$ that hosts $D_{\text{out}}$. More precise idea about the desired maneuver will be given and justified in the next section based on our main assumptions about the surface. In turns, they are slightly enhanced conditions necessary for a robot with a limited turning radius to be capable of coping with contortion of the scanned surface $B$. To come into details of these conditions and assumptions, we need the following notations.
\begin{itemize}
\item $\ts_{b}$ --- the plane tangent to $B$ at point $b \in B$;
\item
$\theta_b \in [0,\pi]$ --- the angle between the normal $N(b)$ and $\vec{h}$;
\item
$S_{b}(V) = - D_{V}N \in \ts_{\bldr}$ --- the shape operator, i.e., minus the derivative of $N$ in the tangent direction $V$;
\item
$\ff{b}{V}{W} = \spr{S_{b}(V)}{W}$ --- the second fundamental form;
\item
$\gamma(h_\ast)= \{\bldr \in B: h(\bldr) = h_\ast\}$ -- the horizontal slice of $B$ at the altitude $h_\ast$; $\gamma_b :=\gamma[h(b)]$;
 \item $\vec{a}_1 \times \vec{a}_2$ --- the cross-product of two vectors $\vec{a}_i \in \br^3$;
\item
$T_0(b)=\frac{N_b \times \vec{h}}{\sin \theta_b}$ --- the unit vector tangential to the slice $\gamma_b$ at $b \in B$ and directed so that $D$ is to the left;
\item
$T_{\bot}(b) := N_b \times T_0(b) \in \ts_b$ --- the second component of orthogonal basis $(T_0,T_{\bot})$ of the tangent plane $\ts_b$;
\item
$\vec{n}_b := \frac{N_b - \cos \theta_b \vec{h}}{\sin \theta_b}$ --- the inner unit normal to the slice $\gamma_b$ in the host horizontal plane;
\item
$b_{\bldr}$ --- the point of $B$ horizontally nearest to $\bldr$.
\end{itemize}
Within the horizontal layer $h(b) \in [h_-,h_+]$, Assumption~\ref{ass.smooth} ensures that $\theta_b \in(0, \pi)$ and $T_0(b)$ is well defined.
\par
As was remarked, we are interested in circumnavigation of $B$ with a small increment $\Delta h$ in altitude between the rounds of circumnavigation. Looking for necessary conditions, we focus on the idealized limit case where $\Delta h =0$ and so the ideal maneuver comes to planar horizontal circumnavigation of a horizontal slice of $B$ with a constant margin of $d_0$. By Lemma~1 in \cite{MTS11} and the associated discussion, the following property is nearly necessary for this maneuver to be feasible at any altitude $h \in [h_-,h_+]$.
\begin{assumption}
\label{ass.mnp}
The ``minimum distance'' point $b_{\bldr} \in B$ is unique for any location $\bldr \in D_{\text{out}}$ at a horizontal distance of $d(\bldr)=d_0$ from $B$, and such locations do exist.
\end{assumption}
This assumption is necessarily fulfilled if any horizontal slice of the domain $D$ is convex.
\par
 The next condition is that the robot's minimum turning radius $R$ should be small enough as compared with contortion of $B$ evaluated by the second fundamental form.
\begin{lemma}
\label{lem.nec1}
 Let the robot can circumnavigate the horizontal slice $\gamma$ of $B$ at any altitude $h\in[h_-,h_+]$, moving in $D_{\text{out}}$ and in the associated horizontal plane, with always maintaining the distance $d_0$ to $\gamma$. Then
 \begin{equation}
 \label{nn.cc}
 R=\frac{v}{\ov{u}} \leq \frac{\sin \theta_b}{|\ff{b}{T_0}{T_0}|} + d_0 \sgn \ff{b}{T_0}{T_0}
 \end{equation}
for any $b \in B$, where $T_0:=T_0(b)$, $a/0:= \infty \; \forall a >0$, and $\sgn 0 := 0$.
\end{lemma}
\par
The proof of this claim employs several technical facts concentrated in the next lemma. To state it, we
recall that $D_V$ is the derivative in direction $V$ and $[\vec{a},\vec{b}, \vec{c}\,] = \sprs{\vec{a}}{\vec{b} \times \vec{c}}$ is the scalar triple product.
\begin{lemma}
At any point $b \in B$ and for any $V \in \ts_b$,
\begin{gather}
\label{tr1}
\sprs{N}{\vec{h}} = \cos \theta, \spr{N}{T_0} = \spr{N}{T_\bot} =0,
\\
\label{tr2}
\sprs{\vec{h}}{T_0} = \sprs{\vec{h}}{\vec{n}} =0, \sprs{\vec{h}}{T_\bot} = - \sprs{N}{\vec{n}} =  - \sin \theta,
\\
\label{tr3}
\sprs{\vec{n}}{T_0} = \sprs{T_\bot}{T_0}=0, \quad \sprs{\vec{n}}{T_\bot} = \cos \theta,
\\
\label{tr4}
D_V \theta = -\ff{}{V}{T_\bot},
\\
\label{tr5}
 D_V \vec{n} = - \frac{\ff{}{V}{T_0}}{\sin \theta}T_0,
\quad
D_V T_0 = \frac{\ff{}{V}{T_0}}{\sin \theta} \vec{n}
\\
\label{tr6}
D_v T_\bot = \ff{}{V}{T_\bot} N - \cot \theta \cdot \ff{}{V}{T_0} T_0.
\end{gather}
\end{lemma}
\pf Formulas \eqref{tr1}---\eqref{tr3} are immediate from the definitions of $\theta, T_0,T_\bot$, and $\vec{n}$ since
\begin{gather}
\nonumber
\sprs{N}{\vec{n}} = \spr{N}{ \frac{N - \cos \theta \vec{h}}{\sin \theta}} = \frac{1 - \cos \theta \sprs{N}{\vec{h}}}{\sin \theta}
 = \frac{1-\cos^2 \theta}{\sin \theta} = \sin \theta,
\\
\nonumber
\sprs{\vec{h}}{\vec{n}} = \spr{\vec{h}}{\frac{N - \cos \theta \vec{h}}{\sin \theta}} = \frac{\sprs{\vec{h}}{N} - \cos \theta}{\sin \theta} =0,
\\
\sprs{\vec{h}}{T_\bot} = \spr{\vec{h}}{N \times T_0} = \left[ \vec{h}, N, T_0 \right] \overset{\text{(a)}}{=} - \left[ T_0, N, \vec{h}\right]
\label{trpr.equal}
= -\big\langle T_0; \underbrace{N \times \vec{h}}_{= \sin \theta T_0} \big\rangle = - \sin \theta;
\\
\nonumber
\sprs{\vec{n}}{T_\bot} = \spr{\frac{N - \cos \theta \vec{h}}{\sin \theta}}{N \times T_0} = -\frac{\cos \theta}{\sin \theta}\spr{\vec{h}}{N \times T_0}
\overset{\text{\eqref{trpr.equal}}}{=\!=} \cos \theta.
\end{gather}
Here (a) holds since the triple scalar product is anti-symmetric. By differentiating $\cos \theta = \sprs{N}{\vec{h}}$, we arrive at
$-D_V \theta \sin \theta = - \spr{S(V)}{\vec{h}} \overset{\text{(b)}}{=} - \spr{S(V)}{\vec{\pi}}$, where $\vec{\pi} := \vec{h} - \cos \theta N$ and (b) holds since $S(V) \in \ts$ and $N$ is normal to the tangent plane $\ts$. Since $\sprs{\vec{\pi}}{N} =0 \Rightarrow  \vec{\pi} \in \ts$, we have
$\vec{\pi} = c T_0 + e T_\bot, c,e \in \br$. Here $c = \sprs{\vec{h} - \cos \theta N}{T_0} = 0$ and $e = \sprs{\vec{h} - \cos \theta N}{T_\bot} = - \sin \theta$  by \eqref{tr1} and \eqref{tr2}. Thus we arrive at \eqref{tr4}.
The first part of \eqref{tr5} holds since $\sprs{\vec{n}}{\vec{h}}=0\;\forall b \Rightarrow \spr{D_V \vec{n}}{\vec{h}} =0$ and $\|\vec{n}\| =1$ $ \forall b  \Rightarrow \spr{D_V \vec{n}}{\vec{n}} =0$ imply that $D_V \vec{n}$ is co-linear with $ \vec{n} \times \vec{h} = \frac{N - \cos \theta \vec{h}}{\sin \theta} \times \vec{h} = \frac{N \times \vec{h}}{\sin \theta} = T_0$, i.e., $D_v \vec{n} = \zeta T_0$, where
\begin{multline*}
\zeta = \spr{D_V \vec{n}}{T_0} = \spr{D_V \frac{N -   \vec{h}\cos \theta}{\sin \theta}}{T_0}
= \spr{\frac{- S(V) + \vec{h} D_V \theta \sin \theta }{\sin \theta}}{T_0} - \spr{\vec{n}}{T_0}D_V \theta \cot \theta
\\
\overset{\text{\eqref{tr2}, \eqref{tr3}}}{=\!=\!=\!=\!=\!=} \spr{\frac{- S(V)}{\sin \theta}}{T_0} .
\end{multline*}
Similar arguments assure that $D_V T_0 = \xi \vec{n}$. Here
\begin{multline*}
\xi = \spr{D_V T_0}{\vec{n}} = \spr{D_V \frac{N \times \vec{h}}{\sin \theta}}{\vec{n}} = - \spr{\frac{S(V) \times \vec{h}}{\sin \theta}}{\vec{n}}
- \spr{T_0}{\vec{n}} \frac{D_V \theta \cos \theta}{\sin \theta} \overset{\text{\eqref{tr3}}}{=} -
\frac{[\vec{n}, S(V), \vec{h}] }{\sin \theta} \overset{\text{(a)}}{=} \frac{[S(V), \vec{n}, \vec{h}] }{\sin \theta}
\\
= \spr{\frac{S(V)}{\sin \theta}}{ \frac{N - \cos \theta \vec{h}}{\sin \theta}\times \vec{h}} = \spr{\frac{S(V)}{\sin \theta}}{ \frac{N \times \vec{h}}{\sin \theta}},
\end{multline*}
where $\frac{N \times \vec{h}}{\sin \theta} = T_0$, which completes the proof of \eqref{tr5}. Finally, $D_VT_\bot$ lies in the plane normal to $T_\bot$ since $\|T_\bot\|=1 \; \forall b$. By \eqref{tr1} and \eqref{tr3}, $T_0$ and $N$ constitute a basis of this plane, and so $D_V T_\bot = \varsigma T_0 + \beta N$. It remains to note that
\begin{gather*}
D_v T_\bot= D_V \left( N \times T_0 \right) = \left( D_V  N \right) \times T_0+ N \times D_V T_0 ,
\\
\varsigma = \spr{D_V T_\bot}{T_0}
= \frac{\ff{}{V}{T_0}}{\sin \theta} \spr{T_0}{N \times \left( \frac{N - \cos \theta \vec{h}}{\sin \theta}\right)}
\\
= -\ff{}{V}{T_0} \cot \theta  \spr{T_0}{\frac{N \times \vec{h}}{\sin \theta}} = -\ff{}{V}{T_0} \cot \theta,
\\
\beta = \spr{D_v T_\bot}{N} = \spr{\left( D_V  N \right) \times T_0}{N} =  - \spr{S(V) \times T_0}{N}
\\
= - \left[ N, S(V), T_0\right] \overset{\text{(a)}}{=} \left[ S(V), N, T_0\right] = \sprs{S(V)}{N \times T_0}
= \sprs{S(V)}{T_\bot} = \ff{}{V}{T_\bot}. \text{\epf}
\end{gather*}
\indent
{\bf Proof of Lemma~\ref{lem.nec1}.} Let $\bro(\sigma)$ be a regular parametric representation of $\gamma$ near $b_{\bldr}$. Here $\sigma$ is the arclength, $\sigma =0$ represents $b_{\bldr}$, and $\sigma$ ascends in the direction of $T_0$. Then $f(\sigma) := 1/2 \|\bro(\sigma) - \bldr\|^2 \geq f(0)$. Hence $\frac{d f(\sigma)}{d \sigma}(0) = \spr{T_0(b_{\bldr})}{b_{\bldr} - \bldr} =0 \Rightarrow b_{\bldr} - \bldr = d_0 \vec{n}_{b_{\bldr}}$. Furthermore,
\begin{gather*}
0 \leq \frac{d^2 f(\sigma)}{d \sigma^2}(0) = \frac{d}{d \sigma} \spr{T_0[\bro(\sigma)]}{ \bro(\sigma)-\bldr}\big|_{\sigma=0}
\overset{\text{\eqref{tr5}}}{=\!=} \spr{\frac{\ff{b_{\bldr}}{T_0(b_{\bldr})}{T_0(b_{\bldr})}}{\sin \theta}\vec{n}_{b_{\bldr}}}{ d_0 \vec{n}_{b_{\bldr}}} + \spr{T_0(b_{\bldr})}{T_0(b_{\bldr})}
\\
= 1+ d_0 \frac{\ff{b_{\bldr}}{T_0(b_{\bldr})}{T_0(b_{\bldr})}}{\sin \theta} =: C(\bldr).
\end{gather*}
Meanwhile, the robot's motion can be parametrized as follows $\bldr(t) = \bro[\sigma(t)] - d_0 \vec{n}_{\rho[\sigma(t)]}$, where $\sigma(t)$ is the parameter of $b_{\bldr(t)}$. With regard to \eqref{1b} and \eqref{tr5}, we have
$$
v \imm = \dot{\bldr} = \dot{\sigma} T_0 + d_0 \dot{\sigma} \frac{\ff{b_{\bldr}}{T_0(b_{\bldr})}{T_0(b_{\bldr})}}{\sin \theta} T_0 = C \dot{\sigma} T_0.
$$
So $\dot{\sigma} = v/C$ and the robot's acceleration $\vec{u}$ is normal to $T_0$ by \eqref{1b}; it also horizontal since
the robot moves in the horizontal plane. It follows that $\vec{u} = u \vec{n}$, where
$$
u = \spr{\vec{u}}{\vec{n}} = \spr{\dot{\imm}}{\vec{n}}\overset{\text{\eqref{tr5}}}{=\!=} \frac{1}{v} \spr{\frac{d(C \dot{\sigma})}{dt} T_0}{\vec{n}}
+ \frac{C \dot{\sigma}^2}{v} \spr{\frac{\ff{}{T_0}{T_0}}{\sin \theta} \vec{n}}{\vec{n}} \overset{\text{\eqref{tr3}}}{=\!=} \frac{v}{C} \frac{\ff{}{T_0}{T_0}}{\sin \theta}.
$$
By putting this into the last relation $|u| = \|\vec{u}\| \leq \ov{u}$ from \eqref{1b} and invoking that $C \geq 0$ by the foregoing, we arrive at \eqref{nn.cc} via elementary transformations. \epf
\par
Thus \eqref{nn.cc} is needed for capacity to maintain $d \equiv d_0$, which at least means that whenever $d=d_0$ and $\dot{d}=0$, a feasible control exists such that $\ddot{d}=0$. It does not come as a surprise that putting $<$ in place of $\leq$ in \eqref{nn.cc} implies the possibility to freely manipulate the sign of $\ddot{d}$ by picking proper feasible controls, thus making the output $d$ locally controllable. In problems of regulating an output to a desired value, the output controllability is conventionally needed at all points on the transient portion of the trajectory. In view of this, we enhance \eqref{nn.cc} by the replacement ``$\leq$'' $\mapsto$ ``$<$'' and extend the resultant inequality on the entire operational zone $D_{\text{op}}$.
\par
For convenience's sake, we define this zone in terms of $h$ and $d$. Since upward and downward scans should be terminated and commenced, respectively, at the altitude $h = h_+$ and the vertical velocity cannot be instantly reversed, transfer from the former to the latter inevitably involves a maneuver at altitudes $h > h_+$. Similarly, the opposite transfer is through altitudes $h < h_-$. So we expand the altitude range from $[h_-,h_+]$ to $[h_- - \Delta_h, h_+ + \Delta_h]$ (with a certain $\Delta_h>0$) in the definition of the operational zone. It is also characterized by the extreme values $d_- < d_+$ taken by $d$ in it, where $\dsafe \leq d_- < d_0 < d_+$. Overall,
$$
D_{\text{op}}:= \{\bldr \in D_{\text{out}}: h_--\Delta_h \leq h(\bldr) \leq h_++\Delta_h, d_- \leq d(\bldr) \leq d_+\},
$$
and  we arrive at the following assumption, which covers Assumption~\ref{ass.mnp} and the second part of Assumptions~\ref{ass.smooth}.
\begin{assumption}
\label{as.cum}
Assumption~{\rm \ref{ass.mnp}} is valid with $h_+:= h_++\Delta_h$ and $h_-:= h_- - \Delta_h$. For any point $\bldr \in D_{\text{op}}$, \eqref{nn.cc} holds with $<$ put in place of $\leq$ and $d_0:= d(\bldr)$.
\end{assumption}
From now on, $B$ denotes the extension of the original $B$ to the lateral piece of the boundary of $D_{\text{op}}$.
\par
Our last assumption addresses the initial location $\bldr_{\text{in}}$ and orientation $\imm_{\text{in}}$ of the robot. To simplify the formulations, we assume that the latter is horizontal $\sprs{\imm_{\text{in}}}{\vec{h}}=0$, which can be always achieved by a preliminary maneuver. Then $\sprs{\vec{u}}{\vec{h}} \equiv 0$ guarantees that the robot remains in the initial horizontal plane. There are two options of such a motion with a maximum and constant actuation
\begin{equation}
\label{max.act}
\|\vec{u}\| \equiv \ov{u},
 \end{equation}
i.e., the left and right turn, respectively. The respective circular paths $C^\pm_{\text{in}}$ are called the {\it initial circles}; the disks $\mathscr{D}^\pm_{\text{in}}$ bounded by them are called the {\it initial disks}. Under the control law proposed in Section~\ref{sec.claw}, the robot initially traces a part of $C^\pm_{\text{in}}$. The last assumption guarantees that the robot does not leave the operation zone when doing so.
\begin{assumption}
\label{ass.inc}
The both initial disks lie in the interior of the operational zone.
\end{assumption}
\section{Specification of the Problem Setup}
\label{sec.specif1}
\begin{define}
\label{def.ffccss}
A maneuver of the robot in $D_{\text{op}}$ is called a {\em full scan of $B$ with preciseness} $\ve_\ast >0$ if for any point $b \in B$, there exists a point $\bldr$ on the robot's path such that $\|b-b_{\bldr}\|< \ve_\ast$.
\end{define}
In fact, the core objective is to perform such a scan with small enough $\ve_\ast$. This leaves much freedom with respect to the paths of both the robot and its projection $b_{\bldr}$. The next definition specifies them by choosing a spiral-like pattern of motion along the surface of interest.
\begin{define}
\label{def.spisc}
A maneuver in $D_{\text{out}}$ is called the {\em pure spiralling scan of $B$ with vertical rate} $\eta \in (-v,v)$ {\em and within the altitude range $[h_-,h_+]$} if the altitude $h$ runs the entire length of this range with the constant rate $\dot{h} = \eta$, and the horizontal distance $d$ from the robot to $B$ is kept equal to $d_0$.
\par
For $\ve,  \ve^d>0$, a maneuver is called an $(\ve,\ve^d)$-{\em approximate spiralling scan of $B$ with vertical rate} $\eta$ {\em and within the altitude range $[h_-,h_+]$} if the above claims hold with $\dot{h} = \eta$ and $d \equiv d_0$ replaced by $|\dot{h} - \eta| \leq \ve$ and $|d -d_0|\leq \ve^d$, respectively, and $|\dot{d}| \leq \ve$.
\end{define}
This scan may be upward ($\eta>0$) or downward ($\eta<0$). The forthcoming Lemma~\ref{lem.anal} shows that if $|\eta|, \ve$, and $\ve_a$ are small enough,
approximate spiralling scan of $B$ with vertical rate $\eta$ and within the altitude range $[h_-,h_+]$ guarantees full scan of $B$ in the sense of Definition~\ref{def.ffccss} with high preciseness. Hence such a scan is a fascinating option.
We are interested in repeatedly performing such scans in the sense of the following definition.
\begin{define}
\label{def.scan}
We say that the {\em surface $B$ is repeatedly scanned by the robot} with the vertical speed $\eta_\ast > 0$
if $d(t) \geq d_{\text{\rm safe}} \; \forall t$ and the time interval $[0,\infty)$ can be partitioned into subintervals $0<t_0^- < t^+_0 < t_1^- < t^+_1 < \ldots$ so that $t^\pm_k \to \infty$ as $k \to \infty$, and the following claims hold:
\begin{enumerate}[{\bf i.}]
\item During any time interval $[t^-_k,t^+_k]$, the robot performs an $(\ve_k,\ve^d_k)$-approximately pure scan with the vertical rate $\eta_k=\pm \eta_\ast$, where $\ve_k, \ve^d_k \to 0$ as $k \to \infty$;
\item During any time interval $I_k:=(t^+_k,t^-_{k+1})$, the altitude $h$ is outside the range $[h_-,h_+]$;
\item The intervals $I_k$ are bounded $\sup_{k} (t^-_{k+1}- t^+_k) < \infty$.
\end{enumerate}
\end{define}
It is easy to see that $\eta_{k+1} = - \eta_k \; \forall k$. The intervals $I_k$ accommodate transients required to pass from an upward scan to downward one, and vice versa; the interval $[0,t^-_0]$ accommodates the transient from the initial state. The convergence $\ve_k, \ve_k^d \to 0$ means that the robot's maneuver becomes a nearly pure scan as time progresses.
\par
{\it It is required to design a controller that ensures repeated scans in the sense of Definition}~\ref{def.scan}.
\par
In conclusion of the section, we flesh out the above remark about the relationship between scans in the sense of Definitions~\ref{def.ffccss} and
\ref{def.spisc}, respectively.
\par
To this end, we suppose that Assumptions~\ref{ass.mnp} and \ref{as.cum} hold, and start with introduction of ``pseudo-cylindrical'' coordinates on $B$. The horizontal slice $\gamma(h_-)$ is a connected compact closed smooth curve. So there is an isometrical isomorphism $\varsigma \in S_0^l \xleftrightarrow{j} \gamma(h_-)$ between the circle $S_0^l$ of a proper radius $l$ and $\gamma(h_-)$. The integral curves of the differential equation $\dot{b} = - T_\bot(b)/\sin \theta_b, b \in B$ on the compact manifold $B$ are extensible at least while remaining in the range of altitudes $[h_-,h_+]$ \cite[Ch.~IV]{Lang02}, whereas $\dot{h}=1$ on them. Hence these curves span the range of altitudes $[h_-,h_+]$ and can be parametrized by $h \in [h_-,h_+]$; let $b=\mathfrak{b}(h,\varsigma), h \in [h_-,h_+]$ be the integral curve starting from $\mathfrak{b}(h_-,\varsigma) = j(\varsigma)$.
Smooth dependence of the integral curve on initial data and uniqueness of the solution of the associated Cauchy problem imply that $(h,l) \in [h_-,h_+] \times S_0^l \to \mathfrak{b}(h,\varsigma) \in B$ is the diffeomorphism between the indicated manifolds. Thus $h$ and $\varsigma$ can be viewed as ``pseudo-cylindrical'' coordinates on $B$.
\begin{lemma}
\label{lem.anal}
Suppose that Assumptions~{\rm \ref{ass.mnp}} and {\rm \ref{as.cum}} are valid.
There exist $\eta_0, \ve_0^d, K, M>0$ such that if
$$
0<|\eta|<\eta_0, \quad \ve < |\eta|/2, \quad \ve^d < \ve^d_0,
$$
the following claims hold for any $(\ve,\ve^d)$-approximate spiralling scan (in $D_{\text{op}}$) of $B$ with vertical rate $\eta$ and within the altitude range $[h_-,h_+]$:
\begin{enumerate}[{\bf i)}]
\item This maneuver is a full scan of $B$ with preciseness $K |\eta|$;
\item The coordinate $\varsigma$ of $b_{\bldr}$ completes no less than $N \geq M\frac{h_+-h_-}{|\eta|} -2$ full runs over the circle $S_0^l$.
\end{enumerate}
\end{lemma}
The proof of this lemma is prefaced by several technical facts.
\begin{lemma}
\label{lem.veloc}
While the robot moves in the operational zone,
\begin{gather}
\label{ddot1}
\dot{h} = v \sprs{\imm}{\vec{h}}, \quad \ddot{h} = v \sprs{\vec{u}}{\vec{h}}, \quad \dot{d} = - v \sprs{\imm}{\vec{n}} - \dot{h} \cot \theta,
\\
\nonumber
\ddot{d} = - v \sprs{\vec{u}}{\vec{n}}- \ddot{h} \cot \theta + \lambda^2 \frac{\ff{}{T_0}{T_0}}{\sin \theta + d \ff{}{T_0}{T_0}}
- 2 \frac{\dot{h} \lambda}{\sin \theta} \frac{ \ff{}{T_\bot}{T_0}}{\sin \theta + d \ff{}{T_0}{T_0}}
\\
+ \frac{\dot{h}^2}{\sin^3 \theta} \left[ \ff{}{T_\bot}{T_\bot} - d \frac{\ff{}{T_0}{T_\bot}^2}{\sin \theta + d \ff{}{T_0}{T_0}}\right],
\label{ddot2}
\end{gather}
where
\begin{equation}
\label{lambda.for}
\lambda := \sgn \sprs{\imm}{T_0} \sqrt{v^2 - \frac{\dot{h}^2}{\sin^2\theta} - 2 \dot{h}\dot{d} \cot \theta - \dot{d}^2}.
\end{equation}
\end{lemma}
\pf The first two formulas in \eqref{ddot1} are immediate from \eqref{1b} since $h = \sprs{\bldr}{\imm}+\text{const}$. The velocity $\vec{\tau}(t) = \frac{d b_{\bldr(t)}}{dt}$ of the point $b_{\bldr(t)} \in B$ nearest to $\bldr(t)$ is tangential to $B$ and so $\vec{\tau} = \zeta T_0 + \xi T_\bot$. By differentiating the equation $\bldr = b - d \vec{n}$ with regard to \eqref{1b} and \eqref{tr5}, we see that
\begin{gather}
\label{vel.vf}
v \imm =  \xi T_\bot - \dot{d} \vec{n} + \underbrace{\left[ \zeta + d \frac{\ff{}{\vec{\tau}}{T_0}}{\sin \theta} \right]}_{\lambda} T_0 ;
\\
\label{dot.hh}
\dot{h} \overset{\text{\eqref{ddot1}}}{=\!=} v \sprs{\imm}{\vec{h}} \overset{\text{\eqref{tr2}}}{=\!=} - \xi\sin \theta \Rightarrow \xi = - \frac{\dot{h}}{\sin \theta};
\\
\nonumber
\dot{d} \overset{\text{\eqref{tr3}}}{=} - v \sprs{\imm}{\vec{n}} + \xi \cos \theta = - v \sprs{\imm}{\vec{n}} - \dot{h} \cot \theta \Rightarrow \text{ \eqref{ddot1}},
\\
\nonumber
v^2 = \sprs{v\imm}{v \imm} \overset{\text{\eqref{tr3}}}{=\!=} \frac{\dot{h}^2}{\sin^2\theta} + 2 \dot{h}\dot{d} \cot \theta + \dot{d}^2 + \lambda^2
\Rightarrow \text{\eqref{lambda.for}},
\\
\nonumber
= \zeta + d \zeta \frac{\ff{}{T_0}{T_0}}{\sin \theta} + d \xi \frac{\ff{}{T_\bot}{T_0}}{\sin \theta}
\Rightarrow
\zeta = \frac{\lambda + d \frac{\dot{h}}{\sin^2 \theta}\ff{}{T_\bot}{T_0} }{1 + d \frac{\ff{}{T_0}{T_0}}{\sin \theta}}.
\label{for.zeta}
\end{gather}
By differentiating \eqref{vel.vf}, we get due to  \eqref{1b}, \eqref{tr5}, \eqref{tr6},
\begin{gather*}
v \vec{u} = \dot{\xi} T_\bot + \xi \left\{  \ff{}{\vec{\tau}}{T_\bot} N - \cot \theta \cdot \ff{}{\vec{\tau}}{T_0} T_0 \right\}
- \ddot{d} \vec{n} + d \frac{\ff{}{\vec{\tau}}{T_0}}{\sin \theta}T_0 + \dot{\lambda} T_0 + \lambda \frac{\ff{}{\vec{\tau}}{T_0}}{\sin \theta} \vec{n};
\\
v \sprs{\vec{u}}{\vec{n}} \overset{\text{\eqref{tr2}, \eqref{tr3}}}{=\!=\!=\!=\!=\!=} - \left[ \frac{\ddot{h}}{\sin \theta}- \xi \ff{}{\vec{\tau}}{T_\bot} \cot \theta\right] \cos \theta
- \ddot{d} + \xi \sin \theta \ff{}{\vec{\tau}}{T_\bot}
+ \lambda \frac{\ff{}{\vec{\tau}}{T_0}}{\sin \theta}
\\
= - \ddot{h} \cot \theta - \ddot{d} + \frac{\xi \ff{}{\vec{\tau}}{T_\bot} + \lambda \ff{}{\vec{\tau}}{T_0}}{\sin \theta}.
\end{gather*}
Here with regard to \eqref{for.zeta},
\begin{gather*}
\xi \ff{}{\vec{\tau}}{T_\bot} + \lambda \ff{}{\vec{\tau}}{T_0}
= \xi^2\ff{}{T_\bot}{T_\bot} + \xi \zeta \ff{}{T_0}{T_\bot} + \lambda \zeta \ff{}{T_0}{T_0} + \lambda \xi \ff{}{T_\bot}{T_0}
\\
 = \xi^2\ff{}{T_\bot}{T_\bot} + \lambda \xi \ff{}{T_\bot}{T_0}
  + \xi  \ff{}{T_0}{T_\bot} \left[ \frac{\lambda \sin \theta - d \xi\ff{}{T_\bot}{T_0} }{\sin \theta + d \ff{}{T_0}{T_0}} \right]
  \\
  + \lambda \ff{}{T_0}{T_0} \left[ \frac{\lambda \sin \theta - d \xi\ff{}{T_\bot}{T_0} }{\sin \theta + d \ff{}{T_0}{T_0}} \right]
  \\
  = \xi^2 \left[ \ff{}{T_\bot}{T_\bot} - d \frac{\ff{}{T_0}{T_\bot}^2}{\sin \theta + d \ff{}{T_0}{T_0}}\right]
  + \lambda^2 \sin \theta \frac{\ff{}{T_0}{T_0}}{\sin \theta + d \ff{}{T_0}{T_0}} + 2 \xi \lambda \frac{\sin \theta \ff{}{T_\bot}{T_0}}{\sin \theta + d \ff{}{T_0}{T_0}}
\end{gather*}
Summarizing, we arrive at \eqref{ddot2}. \epf
\begin{lemma}
If the robot is not vertically oriented and lies in $D_{\text{op}}$, the following relations hold:
\begin{multline}
\label{i.mult}
\spr{\imm}{N} = - \dot{d} \frac{\sin \theta}{v},  \sprs{\imm^\bot_h}{\vec{n}} = \frac{\cot \alpha}{v} \left[ \dot{d} +  \cot \theta \dot{h} \right],
 \imm^\bot_\bot = \frac{\vec{h}\times \imm}{\sin \alpha},
\\
  \sprs{\imm^\bot_\bot}{\vec{n}} = \frac{\lambda}{v \sin \alpha}, \; \cos \alpha = \frac{\dot{h}}{v}, \; \sin \alpha = \frac{\sqrt{v^2-\dot{h}^2}}{v}.
\end{multline}
\end{lemma}
\pf
The first formula in \eqref{i.mult} follows from \eqref{vel.vf} and \eqref{tr1}, \eqref{tr2}. The second formula holds since
\begin{equation*}
\sprs{\imm^\bot_h}{\vec{n}} = \spr{\frac{\vec{h} - \cos\alpha \imm}{\sin \alpha}}{\frac{N - \cos \theta \vec{h}}{\sin \theta}}
= \frac{\sprs{\vec{h}}{N} - \cos \alpha \sprs{\imm}{N} - \cos \theta + \cos \alpha \cos \theta \sprs{\imm}{\vec{h}}}{\sin \alpha \sin \theta}
\end{equation*}
The third relation holds since $\imm^\bot_\bot = \imm^\bot_h \times \imm = \frac{\vec{h} - \cos \alpha \cdot \imm}{\sin \alpha}\times \imm= \frac{\vec{h}\times \imm}{\sin \alpha}$. The the fourth formula is true thanks to \eqref{vel.vf} since
\begin{gather*}
\sprs{\imm^\bot_\bot}{\vec{n}} = \spr{\frac{\vec{h}\times \imm}{\sin \alpha}}{\vec{n}} = \frac{[\frac{N - \vec{h}\cos \theta }{\sin \theta},h,\imm]}{\sin \alpha} = \frac{[N,h,\imm]}{\sin \alpha \sin \theta};
\\
\frac{[\imm,N,h]}{\sin \alpha \sin \theta} = \frac{\spr{\imm}{\frac{N\times h}{\sin \theta}}}{\sin \alpha} = \frac{\spr{\imm}{T_0}}{\sin \alpha}.
\end{gather*}
The last two formulas hold since $\dot{h} = v \sprs{\imm}{\vec{h}} = v \cos \alpha$. \epf
\par
{\bf Proof of Lemma~\ref{lem.anal}.}
As was shown in the proof of Lemma~\ref{lem.veloc}, $\frac{d b_{\bldr}}{dt} = \zeta T_0 - \frac{\dot{h}}{\sin \theta}T_\bot$, where $\zeta$ is given by \eqref{lambda.for} and \eqref{for.zeta}. By Assumption~\ref{as.cum}, $\sin \theta + d \ff{}{T_0}{T_0}>0$.
Since $D_{\text{op}}$ is compact, the above strict inequality and $\sin \theta >0$ imply that the concerned quantities are uniformly (over $D_{\text{op}}$) bounded from below by positive constants.
Since $|\dot{h}| \leq |\eta| + |\dot{h} - \eta| \leq 3|\eta|/2, |\dot{d}| \leq \ve < |\eta|/2$ by Definition~\ref{def.scan}, picking $\eta_0$ and $\ve^d_0$ small enough ensures that in \eqref{lambda.for} and \eqref{for.zeta}, $|\lambda| \geq v/2$ and $|\zeta| \geq \zeta_\ast >0$, where the constant $\zeta_\ast$ depends only on $B,d_0$ and $v$. It follows that the coordinate $\varsigma(t)$ of $b_{\bldr(t)}$ obeys the inequality
$|\dot{\varsigma}| \geq k \zeta_\ast$, where $k>0$ is the coefficient of distortion due to the change of the variables. Due to continuity, $\varsigma$ does not change the direction of motion over $S^l_0$. Thus it completes a full run over $S_0^l$ for no less than $l/(k \zeta_\ast)$ time units. Meanwhile $h$ runs the entire length of $[h_-,h_+]$ with the speed $|\dot{h}| \in [|\eta|/2, 3/2 |\eta|]$, which takes a time $T \geq 2(h_+-h_-)/(3|\eta|)$. Hence any point on the cylinder $[h_-,h_+] \times S_0^l$ lies at a distance of no more than $3/2 |\eta|/(k \zeta_\ast) $ from some point of the path traced by the ``pseudo-cylindrical'' coordinates of $b_{\bldr(t)}$. Since the map $\mathfrak{b}(h,\varsigma)$ is Lipschitz continuous \cite[Ch.~IV]{Lang02}, this implies i). To prove ii), it suffices to note that for time $T$, the coordinate $\varsigma$ makes no less than $\lfloor Tk \zeta_\ast/l \rfloor -1 \geq \lfloor 2(h_+-h_-)k \zeta_\ast/(3|\eta| l) \rfloor -1 \geq M \frac{h_+-h_-}{|\eta|} -2$ complete runs over $S_0^l$, where $M:= 2 k \zeta_\ast/(3l)$. \epf
\section{The Proposed Control Law and its Convergence}
\label{sec.claw}
\begin{figure}[h]
\centering
\subfigure[]{
\scalebox{0.2}{\includegraphics{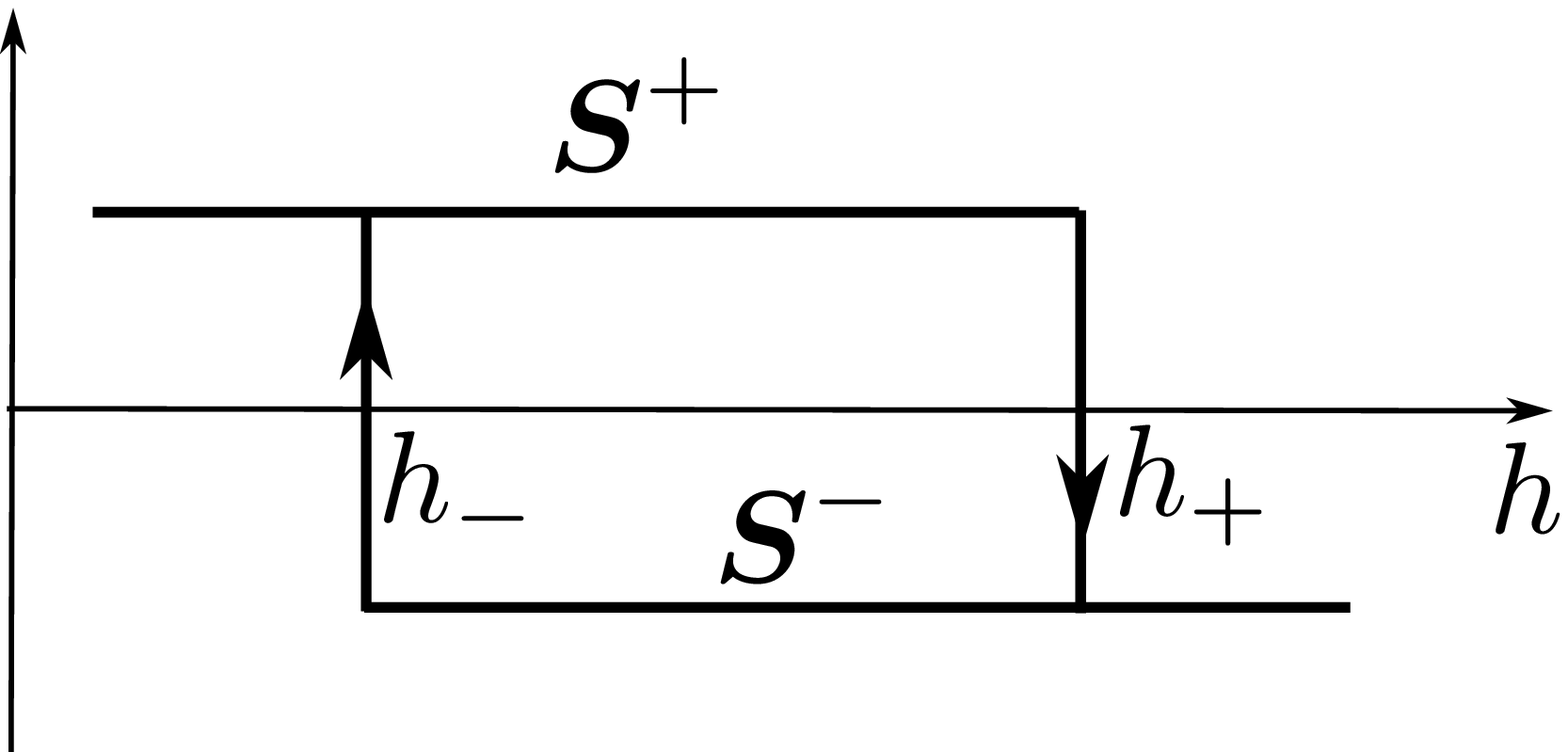}}
}
\subfigure[]{
\scalebox{0.2}{\includegraphics{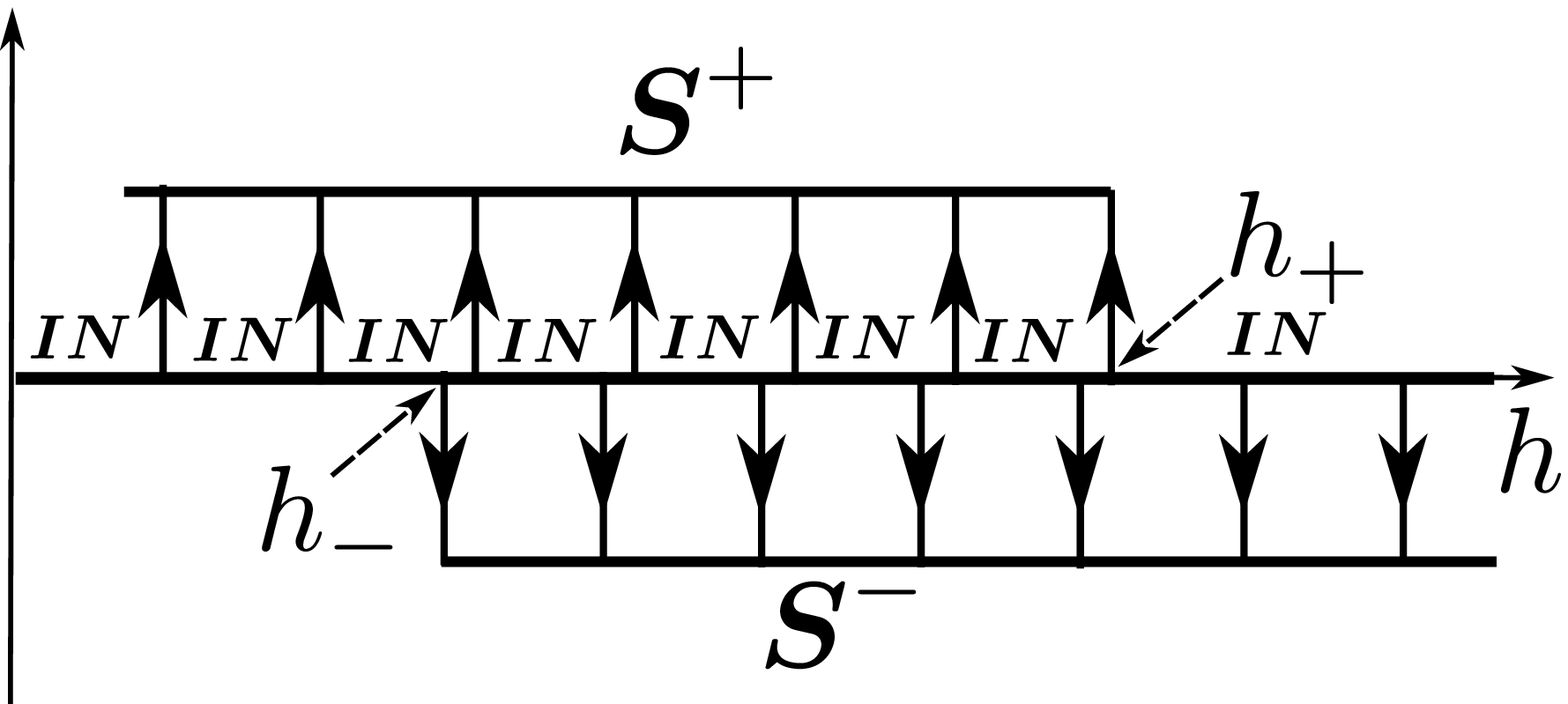}}
}
\caption{Switching logic}
\label{swlogic.fig}
\end{figure}

\setcounter{equation}{0}
The proposed control law assumes that the robot is not vertically oriented, i.e., the angle $\alpha \in [0,\pi]$ between $\imm$ and $\vec{h}$ differs from $0$ and $\pi$. Then the orthogonal projection $\vec{h} - \cos \alpha \imm$ of $\vec{h}$ onto the plane $\imm^\bot$ normal to $\imm$ is nonzero.
Let $\imm^\bot_h := \frac{\vec{h} - \cos \alpha \cdot \imm}{\sin \alpha}$ stand for the normalization of this projection to the unit length, and let $\imm^\bot_\bot := \imm^\bot_h \times \imm$. The pair $(\imm^\bot_h, \imm^\bot_\bot)$ forms an orthonormal basis in $\imm^\bot$ and is computable by the robot in its local reference frame.
To achieve the control objective, we employ a hybrid controller with three discrete states (modes) $\boldsymbol{IN}$, $\boldsymbol{S}^+$, and $\boldsymbol{S}^-$, which mean ``initial maneuver'', ``upward scan'' and ``downward scan'', respectively. The duration $T_{\text{in}}>0$ of the initial maneuver and the scan vertical speed $\eta_\ast >0$ are controller parameters. They are converted into the current vertical rate $\eta$ depending on the mode:
\begin{equation}
\label{def.eta}
\eta := 0 \;\text{in}\; \boldsymbol{IN}, \; \eta := \eta_\ast \;\text{in}\; \boldsymbol{S}^+, \; \eta:= - \eta_\ast \; \text{in}\; \boldsymbol{S^-}.
\end{equation}
The logic of switching the modes is as follows (see Fig.~\ref{swlogic.fig}):
\begin{gather}
\label{switch.s}
\begin{array}{l} \boldsymbol{S}^+ \mapsto \boldsymbol{S}^- \; \text{whenever} \; h > h_+, \\
\boldsymbol{S}^- \mapsto \boldsymbol{S}^+ \; \text{whenever} \; h < h_-.
\end{array}
\\
\boldsymbol{IN} \xmapsto{\text{at time}\; t=T_{\text{in}}}
\begin{cases}
\boldsymbol{S}^+ & \text{if} \; h \leq h_-,
\\
\boldsymbol{S}^\pm & \text{if} \; h \in (h_-,h_+),
\\
\boldsymbol{S}^- & \text{if} \; h \geq h_+.
\end{cases}
\label{intos}
\end{gather}
The initial mode is $\boldsymbol{IN}$; in the middle row from \eqref{intos}, an arbitrary choice between $\boldsymbol{S}^+$ and $\boldsymbol{S}^-$ is performed.
\par
The control law is given by the simple formula
\begin{equation}
\label{c.a}
\vec{u} = - \ov{u}_h \sgn\!\left[ \dot{h} - \eta \right]  \imm^\bot_h
+ \sqrt{\ov{u}^2 - \ov{u}_h^2} \sgn \left[\dot{d} + \chi(d-d_0) \right] \imm^\bot_\bot.
\end{equation}
Here $\chi(\cdot)$ is a linear function with saturation
\begin{equation}
\label{chi}
 \chi(p):=
\begin{cases}
\gamma p & \text{if}\; |p|\leq \delta , \\
\sgn(p)\mu & \text{otherwise},
\end{cases} \qquad \mu :=\gamma \delta,
\end{equation}
and $\ov{u}_h \in (0,\ov{u}), \gamma > 0, \delta >0$ are controller parameters. Since the unit vectors $\imm^\bot_h, \imm^\bot_\bot$, and $\imm$ are mutually perpendicular, the control \eqref{c.a} is feasible, i.e., meets the last two requirements from \eqref{1b}. The time derivatives $\dot{h}$ and $\dot{d}$ can be accessed via, e.g., numerical differentiation.
\par
Solutions of the system closed by the discontinuous controller \eqref{c.a} are meant in the Filippov sense \cite{FIL88}.
\par
Now we show that the control objective can always be achieved by means of the proposed control law under the above nearly unavoidable assumptions.
\begin{theorem}
\label{th.m03} Suppose that Assumptions~{\rm \ref{as.cum}} and {\rm  \ref{ass.inc}} hold. Then there exist controller parameters $\eta_\ast, \ov{u}_h, \mu, \delta$, and $T_{\text{in}}$ such that the following claim holds:
\begin{enumerate}[{\rm (i)}]
\item The proposed controller ensures a repeated scan of the surface $B$ with the requested vertical speed and constantly maintaining a non-vertical orientation $\imm \times \vec{h} \neq 0 \; \forall t$.
\end{enumerate}
 Moreover, for any compact set $Q_{\text{in}} = \{(\bldr, \imm)\}$ of initial states each element of which satisfies Assumption~{\rm  \ref{ass.inc}} and any $\ov{\eta}_\ast>0$, there exist common values of the controller parameters for which {\rm (i)} holds with the vertical speed of scan $\eta_\ast \leq \ov{\eta}_\ast$ whenever the initial state lies in $Q_{\text{in}}$.
\end{theorem}
The proof of the results stated in this section are given in the journal version of the paper submitted to {\it Automatica}. The first claim of Theorem~\ref{th.m03} is a particular case of the second one ($Q_{\text{in}}$ contains only one state).
\par
Now we detail controller tuning under which the conclusion of Theorem~\ref{th.m03} holds.
\par
As a preliminary step, we note that by
Assumption~\ref{ass.inc},
\begin{equation}
\label{strict.in}
\frac{\ov{u}}{v} > \frac{|\ff{}{T_0}{T_0}|}{\sin \theta + d \ff{}{T_0}{T_0}}, \quad \sin \theta + d \ff{}{T_0}{T_0} >0
\end{equation}
in the operational zone.
Since this zone is compact, there exist $k>1$ and $\Delta >0$ such that in this zone,
\begin{equation}
\label{k,delta}
\frac{\ov{u}}{v} > k \frac{|\ff{}{T_0}{T_0}|}{\sin \theta + d \ff{}{T_0}{T_0}} + 2 \Delta.
\end{equation}
After finding such $k$ and $\Delta$, we pick a free design parameter $\varkappa \in (0,1)$.
\par
{\it Tuning $\ov{u}_h  \in (0,\ov{u})$ from} \eqref{c.a}.
To this end, we first modify the definition of the initial disk via putting $\|\vec{u}\| \equiv \sqrt{\ov{u}^2-\ov{u}_h^2}$ in place of \eqref{max.act}. This results in two larger horizontal disks $\mathscr{D}^\pm_{\text{in}}(\ov{u}_h)$ whose boundaries are tangential to $C^\pm_{\text{in}}$ at the initial location of the robot and whose radius $\frac{v}{\sqrt{\ov{u}^2-\ov{u}_h^2}}$ is close to the radius $\frac{v}{\ov{u}}$ of $C^\pm_{\text{in}}$ for $\ov{u}_h \approx 0$. Thanks to Assumption~\ref{ass.inc},
there is $\ov{u}_h^0 \in (0,\ov{u})$ such that the disks $\mathscr{D}^\pm_{\text{in}}(\ov{u}_h)$ lie in the operational zone for all $\ov{u}_h \in (0,\ov{u}_h^0]$. The parameter $\ov{u}_h$ is chosen so that
\begin{gather}
\label{uh.in}
0<\ov{u}_h \leq \min\left\{ \ov{u} \sqrt{1-k^{-1}} ; \ov{u}_h^0 \right\},
\\
\label{uh.in2}
 \ov{u}_h \leq \frac{v \sqrt{1-\kappa}}{\sqrt{k}\left[ \frac{\kappa }{2\sqrt{1-\kappa}}  + | \cot \theta|\right]}\Delta \;\text{in}\; D_{\text{op}}.
\end{gather}
\indent
{\it Tuning $\eta_\ast$ from} \eqref{def.eta} {\it and $\mu, \gamma$ from} \eqref{chi}. These parameters are chosen so that
in the operational zone,
\begin{gather}
\label{tunug}
\frac{\eta_\ast^2}{\sin^2\theta} + 2 \eta_\ast \mu |\cot \theta| + \mu^2 < \kappa v^2,
\\
\label{mn.ch}
2 \eta_\ast \mu \leq (k-1) v^2 |\tan \theta|,
\\
\frac{\eta_\ast^2}{\sin^3 \theta} \left| \ff{}{T_\bot}{T_\bot} - d \frac{\ff{}{T_0}{T_\bot}^2}{\sin \theta + d \ff{}{T_0}{T_0}}\right|  + \gamma \mu
+\frac{  2 v \sqrt{k} }{\sin \theta} \frac{ |\ff{}{T_\bot}{T_0}|}{\sin \theta + d \ff{}{T_0}{T_0}} \eta_\ast
 \leq \frac{v^2 \sqrt{1-\kappa}}{\sqrt{k}}\Delta,
 \label{dec.munu}
\\
\label{eta.deltah}
\eta_\ast < v \sqrt{1 - \left[ 1 - \frac{\Delta_h \ov{u}_h}{v}\right]_+^2}, \;  [a]_+:=\max \{0,a\}.
\end{gather}
Finally, the duration of the initial mode should be chosen so that
\begin{equation}
\label{im.choice}
T_{\text{in}} > \frac{2 \pi}{\sqrt{\ov{u}^2-\ov{u}_h^2}} .
\end{equation}
\par
The recommended choice is feasible since \eqref{uh.in}---\eqref{dec.munu} are satisfied by all small enough $\ov{u}_h, \eta_\ast,\mu$, and $\gamma$.
This can be viewed as a guideline for experimentally tuning the controller. If an a priory knowledge about the domain is available, \eqref{uh.in}---\eqref{dec.munu} can be used for analytically tuning the controller, which is based on estimation of the involved parameters $\theta$ and $\ff{}{\cdot}{\cdot}$ of the surface. In the next section, we illustrate this with simple examples.
\section{Illustrative Examples: Inner and Outer Scans of an Unknown Surface of Revolution}
\label{sec.exmpl}
\setcounter{equation}{0}
Let the surface of interest $B$ be generated by rotating a $C^3$-smooth planar regular curve $C$ around a vertical axis $A$. Also, let the generatrix $C$ lie in a plane containing $A$, do not intersect $A$, and let its tangent $T_C$ be never normal to $A$. It follows that $C$ can be viewed as the graph of the function $\rho(h)$ that gives the distance between $C$ and $A$ at the altitude $h \in I_h:=[h_--\Delta_h, h_+ + \Delta_h]$.
\par
Assuming that the $z$-axis of the world frame is aligned with $A$, we parametrize $x = \rho(h) \cos \phi, y = \rho(h) \sin \phi, z =h$ the surface $B$ by $\phi$ and $h\in I_h$. By using this, it is easy to show (see \cite[Sec.~5.10]{LeClo03} for details) that
\begin{gather*}
N= - \frac{(-1)^\sigma}{\sqrt{1+(\rho^\prime)^2}} \left[ \begin{array}{c} \cos \phi \\ \sin \phi \\ -\rho^\prime\end{array}\right],
\; \theta = \frac{\pi}{2} - (-1)^\sigma \beta,
\\
 \ff{}{T_0}{T_0} = \frac{(-1)^\sigma}{\rho \sqrt{1+(\rho^\prime)^2}} = (-1)^\sigma\frac{\cos \beta}{\rho}, \;  \ff{}{T_0}{T_\bot} = 0,
 \\
 \ff{}{T_\bot}{T_\bot} =  \frac{(-1)^{\sigma+1} \rho^{\prime\prime}}{(1+(\rho^\prime)^2)^{3/2}} = (-1)^{\sigma+1} c.
\end{gather*}
Here $\beta(h)=\arctan \rho^\prime(h)$ is the angle from $\vec{h}$ to $T_C$,
$\sigma =0$ and $\sigma=1$ in the case of outer and inner scan, respectively, and
$c=c(h)$ is the signed curvature of $C$ at the altitude $h$ ($c>0$ at concavities of the graph).
After simple computations, the necessary conditions from Section~\ref{sec.specif} take the form
$R \leq  \rho(h) + (-1)^\sigma d_0$ with an apparent sense: The minimal turning radius $R$ is so small that the robot can horizontally rotate around $A$ at a distance of $\rho + (-1)^\sigma d_0$, which is equivalent to
horizontally following the surface $B$ with the margin $d=d_0$.
\par
Now we comment on analytically tuning the controller, assuming that the available knowledge about the unknown curve $C$ comes to the following estimates:
\begin{itemize}
\item $\ov{\beta}\in [0,\pi/2)$ --- an upper bound on the tangent angle $|\beta|$;
\item $0<\rho_- < \rho_+$ --- a lower and upper estimate of the distance $\rho(h), h\in I_h$ from $A$ to the curve $C$, respectively;
\item $c_+$ --- an upper estimate of the curvature $|c|$ of $C$.
\end{itemize}
Based on these data, the controller should be tuned so that from any initial location $\bldr_{\text{in}}$ with $h(\bldr_{\text{in}}) \in I_h$ and $\rho(\bldr_{\text{in}}) \in [\rho^-_{\text{in}}, \rho^+_{\text{in}}]$, the robot is driven to a horizontal distance of $d_0$ to $B$ and then repeatedly scan $B$ with $d \approx d_0$. Here and in the sequel, $\rho(\bldr)$ is the distance from $\bldr$ to $A$, and $0< \rho^-_{\text{in}} < \rho^+_{\text{in}}$ are given.
\subsection{Outer Scan}
We assume that
\begin{equation}
\rho_{\text{in}}^- > \rho_+ + \dsafe +2R, \quad \rho_{\text{in}}^- > 3R + \rho_+ - \rho_-, \quad \rho_-+d_0 >R.
\label{inner.scan}
\end{equation}
The first two inequalities mean that the initial location $\bldr_{\text{in}}$ is far enough from $A$; the first of them guarantees that $\bldr_{\text{in}}$ is on the correct side of $B$. The last inequality enhances the controllability condition $R \leq  \rho(h) + d_0\; \forall h$.
\par
Now we put
\begin{gather}
\nonumber
d_\ast := \frac{1}{2} \left[ \max \{\dsafe; R - \rho_-\} + \min \{d_0; \rho_{\text{in}}^- - \rho_+  - 2R\}\right],
\quad
\Delta_\rho:= \rho_{\text{in}}^- - \rho_+ - 2R - d_\ast \overset{\text{\eqref{inner.scan}}}{>}0,
\\
\label{uh.def}
u^0_h= \ov{u} \sqrt{1 - \frac{1}{[1+\Delta_\rho/(2R)]^2}}.
\end{gather}
Then we pick free design parameters: $\varkappa \in (0,1)$ and $(k,\Delta)$ from the following triangle in the plane of $(k,\Delta)$'s
\begin{equation}
\label{teke.form}
R^{-1} > \frac{k}{\rho_- + d_\ast } + 2 \Delta, \quad k>1 .
\end{equation}
\begin{proposition}
\label{prop1}
Let the parameters $\ov{u}_h, \eta_\ast, \mu, \gamma >0$ and $T_{\text{\rm in}}$ of the controller be chosen so that \eqref{im.choice} holds and
\begin{gather}
\nonumber
0<\ov{u}_h \leq \min\left\{ \ov{u} \sqrt{1-k^{-1}} ; \ov{u}_h^0; \frac{v \sqrt{1-\kappa}}{\sqrt{k}\left[ \frac{\kappa }{2\sqrt{1-\kappa}}  +  \tan \ov{\beta} \right]}\Delta \right\},
\\
\label{ellips}
\frac{\eta_\ast^2}{\cos^2\ov{\beta}} + 2 \eta_\ast \mu \tan \ov{\beta} + \mu^2 < \kappa v^2,
\\
\label{hyper}
2 \eta_\ast \mu \leq (k-1) v^2 \cot \ov{\beta},
\\
\label{line}
\eta_\ast < v \sqrt{1 - \left[ 1 - \frac{\Delta h \ov{u}_h}{v}\right]_+^2}, \; \frac{\eta_\ast^2  c_+}{\cos^3 \ov{\beta}}
 < \frac{v^2 \sqrt{1-\kappa}}{\sqrt{k}}\Delta,
\\
\label{last}
 \gamma
 \leq \mu^{-1}\left[\frac{v^2 \sqrt{1-\kappa}}{\sqrt{k}}\Delta - \frac{\eta_\ast^2  c_+}{\cos^3 \ov{\beta}}  \right].
\end{gather}
 Then claim {\rm (i)} of Theorem~{\rm \ref{th.m03}} is true whenever $h(\bldr_{\text{\rm in}}) \in I_h, \rho(\bldr_{\text{\rm in}}) \in [\rho^-_{\text{\rm in}}, \rho^+_{\text{\rm in}}]$, and an outer scan is performed.
\end{proposition}
Here the range of $u_h$'s is given in a closed form, whereas \eqref{ellips}, \eqref{hyper}, and \eqref{line} explicitly describe a domain of $(\eta_\ast,\mu)$'s bounded by an ellipse, hyperbola, and a straight line.
After $\mu$ is chosen, \eqref{last} bounds $\gamma$ in an explicit form.
\par
{\bf Proof of Proposition~\ref{prop1}.} By putting $d_- := d_\ast, d_+:=\infty$, it is easy to see that Assumptions~{\rm \ref{as.cum}} and {\rm  \ref{ass.inc}} hold and \eqref{k,delta} takes the form \eqref{teke.form}. It is straightforward to check that the quantity \eqref{uh.def} meets the requirements preceding formula \eqref{uh.in}. By maximizing over $\theta=\pi/2-\beta(h), h\in I_h$, the requirements \eqref{uh.in}---\eqref{mn.ch}, \eqref{eta.deltah} are reduced to those from the proposition's body, except for \eqref{last} and the second part of \eqref{line}. The last inequalities are similarly obtained from \eqref{dec.munu}. \epf
\subsection{Inner Scan} We assume that
\begin{equation}
\rho_- > \rho_{\text{in}}^+ + \dsafe +2R, \quad \rho_{\text{in}}^- > 3R + \rho_+ - \rho_-, \quad \rho_-> d_0 +R.
\label{outer.scan}
\end{equation}
 By the last inequality, the space inside $B$ is ``wide enough''; the first two ones mean that the initial location is far enough from both $B$ and the axis $A$. We put
\begin{gather*}
d_\dagger:= \frac{1}{2} \left[ \max\{d_0; \rho_+-\rho^-_{\text{in}}+ 2R\}+\rho_--R\right],
\\
\Delta_\rho = \min \{ d_\dagger + \rho^-_{\text{in}} - \rho_+ - 2R; \rho_- - \rho_{\text{in}}^+ - \dsafe - 2R\}
\end{gather*}
and define $u^0_h$ by \eqref{uh.def}. The free design parameters: $\varkappa \in (0,1)$ and $(k,\Delta)$ are chosen subject to
\begin{equation}
\label{teke.form1}
R^{-1} > \frac{k}{\rho_- - d_\dagger } + 2 \Delta, \quad k>1 .
\end{equation}
\begin{proposition}
\label{prop2}
Let the parameters $\ov{u}_h, \eta_\ast, \mu, \gamma >0$ and $T_{\text{\rm in}}$
be chosen like in Proposition~{\rm \ref{prop1}}.
 Then claim {\rm (i)} of Theorem~{\rm \ref{th.m03}} is true whenever $h(\bldr_{\text{\rm in}}) \in I_h, \rho(\bldr_{\text{\rm in}}) \in [\rho^-_{\text{\rm in}}, \rho^+_{\text{\rm in}}]$, and an inner scan is performed.
\end{proposition}
\par
The proof is similar to that of Proposition~\ref{prop1}.

 \end{document}